\documentclass{IEEEtran}
\usepackage{cite}
\usepackage{amsmath,amssymb,amsfonts}
\usepackage{multirow}
\usepackage{algorithmic}
\usepackage{amsfonts,cite,times,url,hyperref,bm}
\usepackage{graphicx}
\usepackage{textcomp}
\usepackage{color}
\usepackage{circuitikz}
\usepackage{cancel}

\newtheorem{assumption}{Assumption}

\def\BibTeX{{\rm B\kern-.05em{\sc i\kern-.025em b}\kern-.08em
    T\kern-.1667em\lower.7ex\hbox{E}\kern-.125emX}}
\begin{document}
\title{Distance characteristics for incremental quantities}
\author{Joshua A. Taylor, \IEEEmembership{Senior Member, IEEE}, Alejandro D. Dom\'inguez-Garc\'ia, \IEEEmembership{Fellow, IEEE}
\thanks{This work was supported in part by the National Science Foundation under Grant 2411925. J.~A. Taylor is with the Department of Electrical and Computer Engineering,
       New Jersey Institute of Technology, Newark, NJ, USA. E-mail:
        {\tt\small jat94@njit.edu}. A.~D. Dom\'inguez-Garc\'ia is with the Department of Electrical and Computer Engineering, University of Illinois at 
	Urbana-Champaign, Urbana, IL, USA. E-mail:
        {\tt\small aledan@illinois.edu}}
}

\maketitle

\begin{abstract}
We derive distance relay characteristics in terms of incremental phasors. We use a circuit model of the network to estimate the incremental remote current. If we assume that all sources are stationary, i.e., remain periodic shortly after a fault, then the incremental remote current, and thus the characteristics, do not depend on the real-time voltages or current injections of the sources.
\end{abstract}


\section*{Introduction}

We construct a new type of distance relay characteristic using incremental voltage and current phasors. Incremental quantities are voltages and currents with values from earlier cycles subtracted off, and are used in both the time~\cite{schweitzer2015speed} and phasor domains~\cite{benmouyal2001combined}. In the phasor domain, incremental quantities are sometimes referred to as fault components~\cite{gao2006design,suonan2010novel,ma2011fault,zhang2017integrated}. If the sources remain stationary shortly after a fault, one can show that incremental quantities are driven only by the prefault voltage at the fault point, which can be computed from earlier relay measurements, and not the other sources. When it holds, this assumption makes protection schemes based on incremental quantities independent of the sources.


In this letter, we construct new distance relay characteristics from incremental quantities. As in our earlier work~\cite{taylor2025geometry}, we cast a relay's characteristic as the set of all apparent impedances it can see. The novelty of our work is that we explicitly account for the network structure while preserving computational tractability.

In Section~\ref{sec:network}, we model incremental quantities in terms of the network admittance matrix. These depend on prefault measurements and the network's structure, but not its operating point during the fault. If we assume sources remain periodic for the cycle after a fault, this eliminates source uncertainty and, in turn, remote current uncertainty. The only remaining uncertainty, in the basic distance protection setup, is the fault's location and resistance. In Section~\ref{sec:faultloops}, we derive novel expressions for apparent impedance that depend on incremental quantities. In Section~\ref{sec:discuss}, we interpret our setup and discuss its validity for inverter based resources (IBRs). In particular, stationarity is not a valid assumption for most IBRs today. There is, however, recognition that IBRs will eventually need to respond more consistently to faults for protection systems to operate properly~\cite{schweitzer2026consistency}.

The new characteristics do not have a simple form due to the incremental remote current's nonlinear dependence on the fault location and resistance. In Section~\ref{sec:char}, we construct two novel approximations: a parallelogram based on a point estimate of the remote current, and by simply taking the convex hull of several evaluations of the apparent impedance for different values of the fault location and resistance.

\section{Model}

Let $\delta$ denote the period of sinusoidal voltages and currents. Let $v_{k}(t)\in\mathbb{C}^3$ and $i_{k}(t)\in\mathbb{C}^3$ be the three-phase voltage and current phasors at bus $k$ for the cycle $[t,t+\delta)$. Let $p$ be a positive integer. As in~\cite{schweitzer2015speed}, the corresponding incremental voltage and current are
\[
\tilde{v}_k = v_k(t)-v_k(t-p\delta)\textrm{\quad and\quad}
\tilde{i}_k = i_k(t)-i_k(t-p\delta).
\]
We hereon set $p=1$, and note that all constructions are valid for any positive integer.

Consider a line in the network with positive-sequence series impedance $z$. A relay is at the line's local bus, $\textrm{L}$. The remote bus is denoted $\textrm{R}$. The relay measures the local bus voltages, $v_{\textrm{L}}(t)$, and currents flowing into the line from the local bus, $i_{\textrm{L}}(t)$. $v_{\textrm{R}}(t)$ and $i_{\textrm{R}}(t)$ are the analogous quantities at the remote side, where $i_{\textrm{R}}(t)$ also flows into the line.

There are eleven types of line faults, which we index with the elements of the set $\mathbb{F} = \{\textrm{ag},\textrm{bg},\textrm{cg},
\textrm{ab},\textrm{ac},\textrm{bc},
\textrm{abg},\textrm{acg},\textrm{bcg},
\textrm{abc},\textrm{abcg}
\}.$ For example, ag and ab are phase a-to-ground and phase a-to-phase b faults. Line-to-ground and line-to-line faults are modeled by a resistor to ground and between phases, respectively. As these two faults types are present in the other five faults, in Section~\ref{sec:faultloops} we model only the ag and ab fault loops.

A fault occurs on the line at time $t$. It occurs at normalized location $m_{\textrm{T}}\in\mathbb{R}$, where $m_{\textrm{T}}\in[0,1]$ is on the line and $m_{\textrm{T}}\notin[0,1]$ elsewhere in the network. The fault's resistance is $m_{\textrm{F}}R_{\textrm{F}}$, where $R_{\textrm{F}}>0$ is the maximum possible resistance and $m_{\textrm{F}}\in\left[0,1\right]$. Figure~\ref{fig:circuit} shows the circuit diagram of phases a and b during an ab fault. We write $m$ with no subscript to refer to the pair, $\left(m_{\textrm{T}},m_{\textrm{F}}\right)$. The relay aims to determine if its measurements are consistent with $m_{\textrm{T}}\in[0,1]$, in which case the fault might be on its line. 
\begin{figure}[h]
\begin{circuitikz}
	\draw (0,0) node[above]{$v_{\textrm{L}}^{\textrm{a}}(t)$} to [open, *-*] (8,0) node[above]{$v_{\textrm{R}}^{\textrm{a}}(t)$};
	\draw (0,0) to [generic=$m_{\textrm{T}}z$ , i>_=$i_{\textrm{L}}^{\textrm{a}}(t)$] (4,0)
	to [generic=$(1-m_{\textrm{T}})z$ ,  i_<=$i_{\textrm{R}}^{\textrm{a}}(t)$] (8,0);
	\draw (0,1.5) node[above]{$v_{\textrm{L}}^{\textrm{b}}(t)$} to [open, *-*] (8,1.5) node[above]{$v_{\textrm{R}}^{\textrm{b}}(t)$};
	\draw (0,1.5) to [generic=$m_{\textrm{T}}z$ , i>_=$i_{\textrm{L}}^{\textrm{b}}(t)$] (4,1.5)
	to [generic=$(1-m_{\textrm{T}})z$ ,  i_<=$i_{\textrm{R}}^{\textrm{b}}(t)$] (8,1.5);
\draw (4,0) to [/tikz/circuitikz/bipoles/length=20pt,R=$m_{\textrm{F}}r_{\textrm{F}}$,  i>_=$i_{\textrm{F}}(t)$] (4,1.5)
	to (4,1.5);
	\draw (4,0) to [open, *-*] (4,1.5);
\end{circuitikz}
\caption{Circuit diagram of phases a and b during an ab fault.}
\label{fig:circuit}
\end{figure}


\subsection{Network}\label{sec:network}

When the fault has resistance, the voltage drop depends on the remote current. To model the remote current, we must model the rest of the power network. We split the buses into sets of IBRs, $\mathcal{C}$, SGs, $\mathcal{S}$, junctions and loads, $\mathcal{J}$, and the virtual bus where the fault occurs, $\textrm{F}$. We use subscripts to indicate subvectors; i.e., $v_{\mathcal{C}}(t)$ and $i_{\mathcal{C}}(t)$ are the vectors of IBR voltages and current injections, and $v_{\mathcal{S}}(t)$ and $i_{\mathcal{S}}(t)$ are corresponding SG vectors. We model the SGs as voltage sources and IBRs as current sources in parallel with impedances, i.e., Norton equivalents: $i_{\mathcal{C}}(t)+Y_{\mathcal{C}}v_{\mathcal{C}}(t)$. All junctions satisfy $i_{\mathcal{J}}(t)=-Y_{\mathcal{J}}v_{\mathcal{J}}(t)$, where the admittance is zero for physical junctions. Fault currents satisfy $i_{\textrm{F}}(t)=-Y_{\textrm{F}}^{\eta}v_{\textrm{F}}(t)$, $\eta\in\mathbb{F}$; note that this is only valid for $m_{\textrm{F}}>0$, and that when $m_{\textrm{F}}=0$, the relay's observations do not depend on the network.\footnote{This model breaks down if there is an SG at bus L (or R), in which case $m_{\textrm{T}}=0$ (or $m_{\textrm{T}}=1$) and $m_{\textrm{F}}=0$ corresponds to infinite fault current.}

\begin{assumption}\label{assumption:period}
All sources are periodic over $[t-\delta,t+\delta]$.
\end{assumption}
This means that $v_{\mathcal{S}}(t)=v_{\mathcal{S}}(t-\delta)$ and $i_{\mathcal{C}}(t)=i_{\mathcal{C}}(t-\delta)$. We discuss the validity of this assumption in Section~\ref{sec:discuss}.

Let $Y(m_{\textrm{T}})\in\mathbb{R}^{3(n+1)\times 3(n+1)}$ be the bus admittance matrix, which accounts for the virtual fault bus's location, $m_{\textrm{T}}\in[0,1]$. Substituting for $i_{\textrm{F}}(t)$ and $i_{\mathcal{J}}(t)$, we have
\begin{align}
Y(m_{\textrm{T}}) \begin{bmatrix}
v_{\textrm{F}}(t)\\
v_{\mathcal{J}}(t)\\
v_{\mathcal{C}}(t)\\
v_{\mathcal{S}}(t)
\end{bmatrix}=
\begin{bmatrix}
-Y_{\textrm{F}}^{\eta}(m_{\textrm{F}})v_{\textrm{F}}(t)\\
-Y_{\mathcal{J}}v_{\mathcal{J}}(t)\\
i_{\mathcal{C}}(t)+Y_{\mathcal{C}}v_{\mathcal{C}}(t)\\
i_{\mathcal{S}}(t)
\end{bmatrix}.\label{eq:Ymatt}
\end{align}
Before the fault, we have
\begin{align}
Y(m_{\textrm{T}}) \begin{bmatrix}
v_{\textrm{F}}(t-\delta)\\
v_{\mathcal{J}}(t-\delta)\\
v_{\mathcal{C}}(t-\delta)\\
v_{\mathcal{S}}(t-\delta)
\end{bmatrix}=
\begin{bmatrix}
\bm{0}\\
-Y_{\mathcal{J}}v_{\mathcal{J}}(t-\delta)\\
i_{\mathcal{C}}(t-\delta)+Y_{\mathcal{C}}v_{\mathcal{C}}(t-\delta)\\
i_{\mathcal{S}}(t-\delta)
\end{bmatrix},\label{eq:Ymattdelta}
\end{align}
where $i_{\textrm{F}}(t-\delta)=0$. Subtracting (\ref{eq:Ymattdelta}) from (\ref{eq:Ymatt}), we have
\[
Y(m_{\textrm{T}}) \begin{bmatrix}
\tilde{v}_{\textrm{F}}\\
\tilde{v}_{\mathcal{J}}\\
\tilde{v}_{\mathcal{C}}\\
\bm{0}
\end{bmatrix}=
\begin{bmatrix}
-Y_{\textrm{F}}^{\eta}(m_{\textrm{F}})v_{\textrm{F}}(t)\\
-Y_{\mathcal{J}}\tilde{v}_{\mathcal{J}}\\
Y_{\mathcal{C}}\tilde{v}_{\mathcal{C}}\\
\tilde{i}_{\mathcal{S}}
\end{bmatrix}.
\]
Due to Assumption~\ref{assumption:period}, the source variables have canceled out.

Let $I_{\mathcal{S}}$ be a thin matrix with an identity in the entries corresponding to the buses in $\mathcal{S}$ and zeros elsewhere. We have
\begin{align*}
&\underset{=Y_{\textrm{LHS}}^{\eta}(m)}{\underbrace{\left(Y(m_{\textrm{T}})
+\begin{bmatrix}
Y_{\textrm{F}}^{\eta}(m_{\textrm{F}}) &\bm{0} &\bm{0} & \multirow{4}{*}{$-I_{\mathcal{S}}-Y_{\mathcal{S}}(m_{\textrm{T}})$} \\
\bm{0} & Y_{\mathcal{J}} & \bm{0} & \\
\bm{0} & \bm{0}& -Y_{\mathcal{C}} & \\
\bm{0} & \bm{0} & \bm{0} & 
\end{bmatrix}\right)}}
\begin{bmatrix}
\tilde{v}_{\textrm{F}}\\
\tilde{v}_{\mathcal{J}}\\
\tilde{v}_{\mathcal{C}}\\
\tilde{i}_{\mathcal{S}}
\end{bmatrix}\\
&\quad
=\underset{=Y_{\textrm{RHS}}^{\eta}(m)}{\underbrace{\begin{bmatrix}
-Y_{\textrm{F}}^{\eta}(m_{\textrm{F}})\\
\bm{0}\\
\bm{0}\\
\bm{0}
\end{bmatrix}}}v_{\textrm{F}}(t-\delta).
\end{align*}
Let $\omega^{\eta}(m)=Y_{\textrm{LHS}}^{\eta}(m)^{-1}Y_{\textrm{RHS}}^{\eta}(m)$. Let $D_k$ be such that
\[
v_k(t)=D_k\begin{bmatrix}
v_{\textrm{F}}(t)&
v_{\mathcal{J}}(t)&
v_{\mathcal{C}}(t)&
i_{\mathcal{S}}(t)
\end{bmatrix}^{\top}.
\]
Then $\tilde{v}_k=D_k\omega^{\eta}(m)v_{\textrm{F}}(t-\delta)$ for $k\notin\mathcal{S}$ and $0$ for $k\in\mathcal{S}$. $v_{\textrm{F}}(t-\delta)$ can be written in terms of the relay's observations as $v_{\textrm{F}}(t-\delta) = v_{\textrm{L}}(t-\delta) - m_{\textrm{T}}zi_{\textrm{L}}(t-\delta).$ Let $I_3\in\mathbb{R}^{3\times3}$ denote the identity matrix and
\[
\Omega^{\eta}(m)=\frac{D_{\textrm{R}}-D_{\textrm{F}}}{(1-m_{\textrm{T}})z}\omega^{\eta}(m)\begin{bmatrix}
I_3&-m_{\textrm{T}}zI_3
\end{bmatrix}.
\]
For $m_{\textrm{T}}\in(0,1)$, we can write the incremental remote current (which depends on the fault type, $\eta$) as\footnote{$m_{\textrm{T}}=0$ or $1$ leads to division by zero in this setup. We can instead put the fault resistance at bus L or R, or simply set $m_{\textrm{T}}$ very close to zero or one.}
\begin{align}
\sigma^{\eta}(m) &= \frac{\tilde{v}_{\textrm{R}}-\tilde{v}_{\textrm{F}}}{(1-m_{\textrm{T}})z}\nonumber\\
&=\Omega^{\eta}(m)
\begin{bmatrix}
v_{\textrm{L}}(t-\delta)\\i_{\textrm{L}}(t-\delta)
\end{bmatrix}.\label{eq:tildeieta}
\end{align}

For a given value of $m$, $\sigma^{\eta}(m)$ is easy to compute because $\Omega^{\eta}\left(\hat{m}\right)$ can be computed offline for each $\eta\in\mathbb{F}$, and $v_{\textrm{L}}(t-\delta)$ and $i_{\textrm{L}}(t-\delta)$ are known one or more cycles in advance.

\subsection{Fault loops}\label{sec:faultloops}

For each $\eta\in\mathbb{F}$, the relay computes an apparent voltage, $v_{\textrm{A}}^{\eta}(t)$, and current, $i_{\textrm{A}}^{\eta}(t)$. Let $\psi^{\eta}$ be such that $v_{\textrm{A}}^{\eta}(t)=\psi^{\eta}v_{\textrm{L}}(t)$; we state $v_{\textrm{A}}^{\eta}(t)$, $i_{\textrm{A}}^{\eta}(t)$, and $\psi^{\eta}$ for ag and ab fault loops below. The incremental versions, $\tilde{v}_{\textrm{A}}^{\eta}$ and $\tilde{i}_{\textrm{A}}^{\eta}$, are defined analogously. The apparent impedance and apparent incremental impedance seen by the relay are $z_{\textrm{A}}^{\eta}(m)=v_{\textrm{A}}^{\eta}(t)/i_{\textrm{A}}^{\eta}(t)$ and $\tilde{z}_{\textrm{A}}^{\eta}(m)=\tilde{v}_{\textrm{A}}^{\eta}/\tilde{i}_{\textrm{A}}^{\eta}$.

\subsubsection{Line-to-ground}\label{sec:loop:lg}
Let $z^0$ be the zero-sequence impedance of the line and $k = z^0/z - 1$ the zero-sequence compensation factor. The apparent voltage and current are $v_{\textrm{A}}^{\textrm{ag}}(t) = v_{\textrm{L}}^{\textrm{a}}(t)$ and $i_{\textrm{A}}^{\textrm{ag}}(t) = i_{\textrm{L}}^{\textrm{a}}(t)+ki_{\textrm{L}}^0(t)$,
where $i_{\textrm{L}}^0(t)$ is the zero sequence current. Here $\psi^{\textrm{ag}}=[1,0,0]$.

KVL from the relay to ground gives
\begin{align*}
v_{\textrm{L}}^{\textrm{a}}(t)&=  m_{\textrm{T}} z \left(i_{\textrm{L}}^{\textrm{a}}(t)+ki_{\textrm{L}}^0(t)\right) +  m_{\textrm{F}} r_{\textrm{F}} \left(i_{\textrm{L}}^{\textrm{a}}(t)+ i_{\textrm{R}}^{\textrm{a}}(t)\right)\\
&=m_{\textrm{T}} z \left(i_{\textrm{L}}^{\textrm{a}}(t)+ki_{\textrm{L}}^0(t)\right) +  m_{\textrm{F}} r_{\textrm{F}} \left(\tilde{i}_{\textrm{L}}^{\textrm{a}}+ \psi^{\textrm{ag}}\sigma^{\textrm{ag}}(m)\right)\nonumber\\
&\quad + m_{\textrm{F}} r_{\textrm{F}} \left(i_{\textrm{L}}^{\textrm{a}}(t-\delta)+ i_{\textrm{R}}^{\textrm{a}}(t-\delta)\right)\\
&=m_{\textrm{T}} z \left(i_{\textrm{L}}^{\textrm{a}}(t)+ki_{\textrm{L}}^0(t)\right) +  m_{\textrm{F}} r_{\textrm{F}} \left(\tilde{i}_{\textrm{L}}^{\textrm{a}}+ \psi^{\textrm{ag}}\sigma^{\textrm{ag}}(m)\right),
\end{align*}
where the last line is due to the fact that $i_{\textrm{L}}^{\textrm{a}}(t-\delta)+ i_{\textrm{R}}^{\textrm{a}}(t-\delta)=0$, i.e., the fault current is zero before the fault. Dividing through by $i_{\textrm{A}}^{\textrm{ag}}(t)$, we have
\begin{align}
z_{\textrm{A}}^{\textrm{ag}}(m) &=  m_{\textrm{T}} z  +  m_{\textrm{F}} r_{\textrm{F}} \frac{\tilde{i}_{\textrm{L}}^{\textrm{a}}+ \psi^{\textrm{ag}}\sigma^{\textrm{ag}}(m)}{i_{\textrm{L}}^{\textrm{a}}(t)+ki_{\textrm{L}}^0(t)}.\label{eq:zAag}
\end{align}
A similar derivation leads to the apparent incremental impedance:
\begin{align*}
\tilde{z}_{\textrm{A}}^{\textrm{ag}}(m) &=  m_{\textrm{T}} z  +  m_{\textrm{F}} r_{\textrm{F}} \frac{\tilde{i}_{\textrm{L}}^{\textrm{a}}+ \psi^{\textrm{ag}}\sigma^{\textrm{ag}}(m)}{\tilde{i}_{\textrm{L}}^{\textrm{a}}+k\tilde{i}_{\textrm{L}}^0}.
\end{align*}

\subsubsection{Line-to-line}\label{sec:loop:ll}
The apparent voltage and current are $v_{\textrm{A}}^{\textrm{ab}}(t) = v_{\textrm{L}}^{\textrm{a}}(t)-v_{\textrm{L}}^{\textrm{b}}(t)$ and $i_{\textrm{A}}^{\textrm{ab}}(t) = i_{\textrm{L}}^{\textrm{a}}(t)-i_{\textrm{L}}^{\textrm{b}}(t)$. Here, $\psi^{\textrm{ab}}=[1,-1,0]$. Following the same steps as Section~\ref{sec:loop:lg}, we obtain:
\begin{align}
z_{\textrm{A}}^{\textrm{ab}}(m) &=  m_{\textrm{T}} z  +  \frac{m_{\textrm{F}} r_{\textrm{F}}}{2}\left(  \frac{\tilde{i}_{\textrm{L}}^{\textrm{a}}- \tilde{i}_{\textrm{L}}^{\textrm{b}}+\psi^{\textrm{ab}}\sigma^{\textrm{ab}}(m)}{i_{\textrm{L}}^{\textrm{a}}(t)- i_{\textrm{L}}^{\textrm{b}}(t)}\right)\label{eq:zAab}\\
\tilde{z}_{\textrm{A}}^{\textrm{ab}}(m) &=  m_{\textrm{T}} z  +  \frac{m_{\textrm{F}} r_{\textrm{F}}}{2}\left( 1+ \frac{\psi^{\textrm{ab}}\sigma^{\textrm{ab}}(m)}{\tilde{i}_{\textrm{L}}^{\textrm{a}}- \tilde{i}_{\textrm{L}}^{\textrm{b}}}\right).\nonumber
\end{align}

For both fault loops, the last term of $\tilde{z}_{\textrm{A}}^{\eta}(m)$ could be near $0/0$ when there is no fault. For this reason, we hereon focus on the apparent impedance, $z_{\textrm{A}}^{\eta}(m)$.

\subsection{Interpretation and validity}\label{sec:discuss}

Through (\ref{eq:tildeieta}), the incremental remote current depends on the earlier measurements, $v_{\textrm{L}}(t-\delta)$ and $i_{\textrm{L}}(t-\delta)$, and the network structure. Together with (\ref{eq:zAag}) and (\ref{eq:zAab}), this implies that after the fault, the apparent impedance does not depend on the voltages and currents elsewhere in the network, i.e., its operating point. This parallels the development in~\cite{schweitzer2015speed}, which showed that incremental quantities depend on the prefault voltage at the fault point and Th\'{e}venin impedances of the network.


We now discuss the validity of Assumption~\ref{assumption:period}. It is standard for SGs, whose prefault voltages persist for several cycles after a fault. While there is precedent for modeling IBRs as current sources~\cite{banaiemoqadam2019control} and Norton equivalents~\cite{wieserman2014fault}, Assumption~\ref{assumption:period} is less appropriate for an IBR because its current magnitude can change suddenly after a fault, upon which the inverter will go into current limiting mode. This increase in magnitude, however, is typically smaller than for an SG. An IBRs frequency and phase angle can also change significantly after a fault~\cite{kasztenny2022distance}. This behavior depends on its control algorithm, which can vary case by case. As such, Assumption~\ref{assumption:period} and the implication that $i_{\mathcal{C}}(t)=i_{\mathcal{C}}(t-\delta)$ are generally not realistic today.

That said, it appears that for protection systems to operate effectively in IBR-rich grids, the IBRs should behave more predictably for several cycles after a fault~\cite{schweitzer2026consistency}. In the future, standards enforcing such behavior could make Assumption~\ref{assumption:period} more realistic. If not, we believe that our technical approach could be useful for incorporating other information about IBR fault behavior. For example, we could add uncertainty to $i_{\textrm{C}}(t)$. Such a scheme would still be operating point-independent while allowing for sources to vary after faults. Developing efficiently computable characteristics, e.g., via zonotopes as in~\cite{taylor2025geometry}, is a topic of future work.

\section{Uncertainty and characteristics}\label{sec:char}

Let $\mathcal{Z}^{\eta} = \left\{\left.z_{\textrm{A}}^{\eta}(m) \;\right|\;m_{\textrm{T}}\in[0,1],\; m_{\textrm{F}}\in[0,1]\right\}.$ If a fault of type $\eta$ occurs, then $z_{\textrm{A}}^{\eta}(m)\in\mathcal{Z}^{\eta}$ for any realization of $m\in[0,1]\times[0,1]$. $\mathcal{Z}^{\eta}$ thus corresponds to an overreaching characteristic. Unfortunately, $\mathcal{Z}^{\eta}$ does not appear to have a simple form due to the nonlinear dependency of $\Omega^{\eta}(m)$ on $m$. We instead suggest the following two approximations.

\subsection{Point estimate}\label{sec:point}
Let
\begin{subequations}
\label{eq:zAcert}
\begin{align}
\hat{z}_{\textrm{A}}^{\textrm{ag}}\left(m,\hat{m}\right) &=  m_{\textrm{T}} z  +  m_{\textrm{F}} r_{\textrm{F}} \frac{\tilde{i}_{\textrm{L}}^{\textrm{a}}+ \psi^{\textrm{ag}}\sigma^{\textrm{ag}}\left(\hat{m}\right)}{i_{\textrm{L}}^{\textrm{a}}(t)+ki_{\textrm{L}}^0(t)}\\
\hat{z}_{\textrm{A}}^{\textrm{ab}}\left(m,\hat{m}\right) &=  m_{\textrm{T}} z  +  \frac{m_{\textrm{F}} r_{\textrm{F}}}{2}\left(  \frac{\tilde{i}_{\textrm{L}}^{\textrm{a}}- \tilde{i}_{\textrm{L}}^{\textrm{b}}+\psi^{\textrm{ab}}\sigma^{\textrm{ab}}\left(\hat{m}\right)}{i_{\textrm{L}}^{\textrm{a}}(t)- i_{\textrm{L}}^{\textrm{b}}(t)}\right).
\end{align}
\end{subequations}
We have evaluated the remote current, $\sigma^{\eta}\left(\hat{m}\right)$, at a nominal value of the uncertainty, $\hat{m}$, and left the rest of each expression as they were. Define
\begin{align*}
\hat{\mathcal{Z}}^{\eta}\left(\hat{m}\right) &= \left\{\left.\hat{z}_{\textrm{A}}^{\eta}\left(m,\hat{m}\right) \;\right|\;m_{\textrm{T}}\in[0,1],\; m_{\textrm{F}}\in[0,1]\right\}.
\end{align*}
$\hat{\mathcal{Z}}^{\eta}\left(\hat{m}\right)$ is a parallelogram approximation of $\mathcal{Z}^{\eta}$. It is the Minkowski sum of two line segments, which correspond to the two terms that make up $\hat{z}_{\textrm{A}}^{\eta}\left(m,\hat{m}\right)$ in (\ref{eq:zAcert}). The first line segment is in the direction of $z$. The second line segment depends on $\hat{m}$; a reasonable choice is $\hat{m}_{\textrm{T}}=0.5$ and $\hat{m}_{\textrm{F}}=1$.

\subsection{Convex hull approximation}\label{sec:ch}

Let $\hat{m}_k\in[0,1]\times[0,1]$, $k=1,...,n$. Each $\hat{m}_k$ is a potential realization of $m$, and $z_{\textrm{A}}^{\eta}\left(\hat{m}_k\right)$ is the corresponding apparent impedance for fault type $\eta\in\mathbb{F}$. Let $\textrm{ch}(\cdot)$ denote the convex hull operator. Define
\[
\mathcal{Z}^{\eta}_{\textrm{ch}}=\textrm{ch}\left\{\bigcup_{k=1}^n z_{\textrm{A}}^{\eta}\left(\hat{m}_k\right)
\right\}.
\]
With enough different values of $\hat{m}_k$, $\mathcal{Z}^{\eta}_{\textrm{ch}}$ will roughly contain $\mathcal{Z}^{\eta}$. As each $z_{\textrm{A}}^{\eta}\left(\hat{m}_k\right)$ is a point in the complex plane, taking the convex hull is tractable, e.g., via the gift wrapping algorithm~\cite{jarvis1973identification}. A simple choice for the $\hat{m}_k$ is the four corners, $[0,0],[0,1],[1,0]$,and $[1,1]$.




\section{Example}

\begin{figure}
\includegraphics[width=\columnwidth]{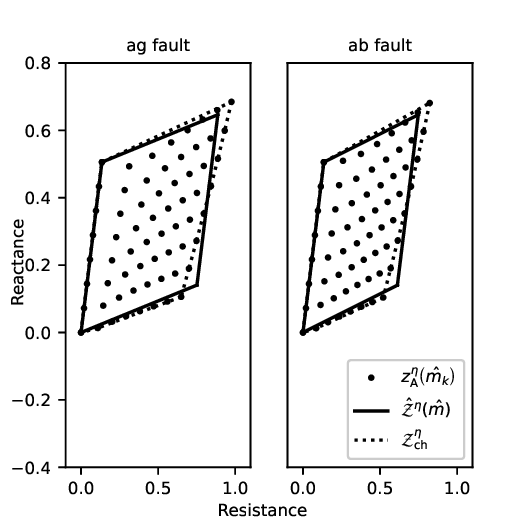}
\caption{Characteristics and impedances for ag and ab faults.}
\label{fig:chars}
\end{figure}

Figure~\ref{fig:chars} shows ag and ab characteristics for the test system in~\cite{taylor2025geometry}. We use $\hat{m}_{\textrm{T}}=0.5$ and $\hat{m}_{\textrm{F}}=1$ for the point estimate in Section~\ref{sec:point}, as well as to generate the measurements $i_{\textrm{L}}(t-\delta)$ and $i_{\textrm{L}}(t)$. For the convex hull approximation in Section~\ref{sec:ch}, we used a $8\times8$ grid over the unit square.


We also show the impedances corresponding to the grid points used to create the convex hull approximation, $z_{\textrm{A}}^{\eta}\left(\hat{m}_k\right)$, which illustrate the exact characteristic, $\mathcal{Z}^{\eta}$. Each $\mathcal{Z}^{\eta}$ has nonconvexity on the top left and very slightly on the bottom right. This means that the convex hull approximation is conservative, i.e., includes unrealizable impedances in these regions. The plots indicate that this error is small.

The convex hull characteristics took roughly 17 milliseconds to compute. This could be sped up substantially by only using $\hat{m}_k$ on the perimeter of $[0,1]\times[0,1]$. For example, using only the four corners takes one millisecond. The point estimates take 0.7 milliseconds to compute and are reasonably close to the convex hull approximations.

\bibliographystyle{IEEEtran}
\bibliography{MainBib,JATBib}

\end{document}